\documentclass[twoside,11pt,titlepage]{amsart}
\usepackage{amsmath} 
\usepackage{amsthm} 
\usepackage{amssymb}	
\usepackage{graphicx} 
\usepackage{multicol} 
\usepackage[marginratio=1:1,height=584pt,width=360pt,tmargin=117pt]{geometry}
\usepackage{hyperref}
\usepackage{tikz}
\usepackage{amscd}

\author{Modjtaba Shokrian Zini, Zhenghan Wang}
\title{Mixed-state TQFTs}

\makeatletter
\newcommand\Author{Modjtaba Shokrian Zini \& Zhenghan Wang}
\let\Title\@title
\def\ps@mystyle{%
      \let\@oddfoot\@empty\let\@evenfoot\@empty
      \def\@evenhead{\makebox[0pt][l]{\thepage}\hfill\Author\hfill}%
      \def\@oddhead{\hfill\Title\hfill\makebox[0pt][l]{\thepage}}%
      \let\@mkboth\markboth}
\makeatother
\makeatletter 
\g@addto@macro{\endabstract}{\@setabstract}
\newcommand{\authorfootnotes}{\renewcommand\thefootnote{\@fnsymbol\c@footnote}}%
\makeatother
\makeatletter
\renewcommand{\maketitle} 
{ \begingroup \vskip 10pt \begin{center} \large {\bf \@title}
	\vskip 10pt \large \@author \hskip 20pt \@date \end{center}
  \vskip 10pt \endgroup \setcounter{footnote}{0} }
\makeatother 

\let\baraccent=\= 
\renewcommand{\=}[1]{\stackrel{#1}{=}} 

\theoremstyle{definition}
\newtheorem{dfn}{Definition}
\theoremstyle{remark}

\usepackage{url,times,wasysym,geometry,indentfirst,rotating,tikz}
\title{\LARGE{\bf
\textsc{Mixed-state TQFTs}}}

\begin{document}
\begin{center}
  \LARGE 
  \maketitle \par \bigskip

  \normalsize
  \authorfootnotes
  Modjtaba Shokrian Zini \footnote{mshokrianzini@pitp.ca}\textsuperscript{,\hyperref[1]{1}}, Zhenghan Wang \footnote{zhenghwa@microsoft.com}\textsuperscript{,\hyperref[2]{2}} \par 
  \bigskip
\end{center}

\begin{abstract}
In this short note, we propose a generalization of Atiyah type TQFTs from pure states to mixed states in the sense that the Hilbert space of pure states associated to a space manifold is replaced by a quantum coherent space related to density matrices. Atiyah type TQFT is a symmetric monoidal functor from the Bord category of manifolds to the category Vec of finite dimensional vector spaces.  In this paper, we define mixed-state TQFTs by replacing the target category Vec by QCS--the category of quantum coherent spaces, then a mixed-state TQFT is simply a symmetric monoidal functor from Bord to QCS.  We also discuss how to construct interesting examples from subsystem quantum error correction codes beyond the trivial ones--all Atiyah type TQFTs. 
\end{abstract}

\section{Introduction}

The pure states of a closed quantum system are modeled by vectors of a Hilbert space, but in practice no quantum systems are really closed (with the possible exception of the universe).  For open quantum systems, their mixed-states are modeled by density matrices, which do not form a vector space.  
Quantum coherent spaces (QCSs) are introduced by Girard to form a new model of linear logic \cite{girard2011blind}, and are closely related to density matrices.   
A QCS is a subset of Hermitian matrices which is equal to its double polar.
QCSs form a symmetric monoidal category \cite{baratella2010quantum}, which is an appealing feature.  In this short note, we propose a generalization of Atiyah type TQFTs from pure states to mixed states in the sense that the Hilbert space of pure states associated to a space manifold is replaced by a quantum coherent space related to density matrices. Atiyah type TQFT is a symmetric monoidal functor from the $\mathbb{B}$ord category of oriented manifolds and their bordisms to the category $\mathbb{V}$ec of finitely dimensional vector spaces and their linear maps \cite{atiyah1988topological}.  In this paper, we define mixed-state TQFTs by replacing the target category $\mathbb{V}$ec by $\mathbb{Q}$cs--the category of quantum coherent spaces, then a mixed-state TQFT is simply a symmetric monoidal functor from $\mathbb{B}$ord to $\mathbb{Q}$cs.  We also discuss how to construct interesting examples from subsystem quantum error correction codes beyond the trivial ones--all Atiyah type TQFTs. 

The Hilbert spaces $V(Y)$ of the Turaev-Viro type TQFTs are quantum error correction codes \cite{qiu2020ground}, and in general we expect many subsystem quantum error correction codes can be made into MS-TQFTs.  We leave the construction of more interesting MS-TQFTs to the future. There was a relevant discussion with our work in the talk on positive TQFTs \cite{Oeckl16}.

All TQFTs in this note will be unitary. 

\section{Category of quantum coherent spaces}

\subsection{Quantum coherent space}

\subsubsection{Polar definition}

Given a finite dimensional Hilbert space $V$ of complex dimension=$n$, the set $H(V)$ of Hermitian matrices form a real Hilbert space of dimension=$n^2$ with the trace norm $<f,g>=tr(fg), f,g\in H(V)$.

\begin{dfn}

\begin{enumerate}
    \item Two Hermitian matrices $f,g\in H(V)$ are polar to each other if $0\leq tr(fg)\leq 1$.
    \item Given a set $C\subseteq H(V)$, the polar of $C$ is the subset $\sim C$ of $H(V)$ such that $\sim C=\{g| g \; \textrm{is polar to every}\; f\in C\}$.
    \item A subset $C\subseteq H(V)$ is a QCS if $\sim \sim C=C$, and we will refer to the Hilbert space $V$ as the carrier of $C$.
\end{enumerate}

\end{dfn}

The canonical example of a QCS is the unit ball in positive Hermitian matrices with trace norm $\mathbb{D}(X)=\{f|f\geq 0, ||f||\leq 1\}$. The polar of  $\mathbb{D}(X)$ is another QCS $\mathbb{P}(X)=\{f|f\geq 0, tr(f)\leq 1\}$.

\subsubsection{Relation to density matrices}

Just as we usually work with the Hilbert spaces of unnormalized pure states, it is easier to work with unnormalized density matrices---positive Hermitian matrices with norm $\leq 1$.

\subsection{The category $\mathbb{Q}$cs}

Given two QCSs $C, D$ with carriers $V$ and $W$, the morphism set from $C$ to $D$ is the set $Hom(C,D)=\{f: f\in H(V\otimes W), f(c)\in D \;\textrm{for every} \; c\in C\}$, where $f$ is naturally identified with a linear map from $H(V)$ to $H(W)$ \cite{baratella2010quantum}.

The tensor product of $C$ and $D$ is $C\otimes D=\sim \sim \{c\otimes d|c\in C,d\in D\}$. The tensor unit is $\bf{1}=\mathbb{R_+}$, the set of $1\times 1$ positive Hermitian matrices.

With this tensor product, QCSs indeed form a symmetric monoidal category as proven in \cite{baratella2010quantum}.
\section{Mixed-state TQFTs}

\subsection{Definition}

\begin{dfn}
An $(n+1)$-Mixed-State TQFT (MS-TQFT) is a symmetric monoidal functor from $\mathbb{B}$ord to $\mathbb{Q}$cs.
\end{dfn}

Obviously, every Atiyah type TQFT gives rise to an MS-TQFT simply by post-composing the map from the Hilbert space $V(Y)$ to the canonical QCS $\mathbb{D}(V(Y))$.

\subsection{Examples of MS-TQFTs}

In real systems, the topological degree of freedom of topological quantum systems would be entangled with other degree of freedom such as environments or measuring devices.  Tracing out all other degree of freedom would result in an MS-TQFT.  Also in real systems, the temperature cannot be absolutely zero as required in ideal topological systems, hence thermal effect could also put the TQFT into a mixed one.  While we will not discuss such realistic mixed-state TQFTs, we point out a non-trivial classe of MS-TQFTs.

\subsubsection{Subsystem of TQFTs}

In $(2+1)$-D Turaev-Viro (TV) type TQFTs, if the input category $\mathcal{B}$ is modular, then the Hilbert space $V_{TV}(Y)$ for each closed surface $Y$ splits into a tensor product $V_{TV}(Y)\cong V_{RT}(Y)\otimes V_{RT}^*(Y)$, where $V_{RT}(Y)$ is the Hilbert space associated to $Y$ from the Reshetikhin-Turaev (RT) TQFT of $\mathcal{B}$. If the subsystem $ V_{RT}^*(Y)$ is traced out, then we obtain an MS-TQFT with $V_{MS-TQFT}(Y)=\mathbb{D}(V_{RT}(Y))$. 

More generally, if an RT TQFT is associated with an MTC  $\mathcal{B}$ which splits into $\mathcal{B}_1\boxtimes \mathcal{B}_2$, then tracing out either system would give rise to an MS-TQFT. A concrete example would be the $SU(2)_k$ Witten-Chern-Simons TQFT for level $k \equiv \pm 1 \pmod 4$.

\vspace{0.5cm}
\noindent \textbf{Acknowledgments.} 
Z.W. is partially supported by NSF grants
FRG-1664351, CCF 2006463, and ARO MURI contract W911NF-20-1-0082.  

\bibliographystyle{unsrt}
\bibliography{references}

\address{\textsuperscript{1\label{1}}Perimeter Institute for Theoretical Physics, Waterloo, Ontario N2L 2Y5, Canada}
\address{\textsuperscript{2\label{2}}Microsoft Station Q and Dept of Mathematics, University of California,
Santa Barbara, CA 93106-6105, U.S.A.} 
\end{document}